\documentclass[11pt]{article}
\usepackage{amsfonts}
\usepackage{amsmath}

\newtheorem{theorem}{Theorem}

\newtheorem{lemma}[theorem]{Lemma}

\newenvironment{proof}[1][Proof]{\noindent\textbf{#1.} }{\ \rule{0.5em}{0.5em}}
\input{tcilatex}
\begin{document}

\title{On the solvability of third-order three point systems of differential
equations with dependence on the first derivative\thanks{%
First author was supported by National Founds through FCT-Funda\c{c}\~{a}o
para a Ci\^{e}ncia e a Tecnologia, project SFRH/BSAB/114246/2016}}
\author{Feliz Minh\'{o}s$^{(\dag )}$ and Robert de Sousa$^{(\ddag )}$ \\
$^{(\dag )}${\small Departamento de Matem\'{a}tica, Escola de Ci\^{e}ncias e
Tecnologia,}\\
{\small Centro de Investiga\c{c}\~{a}o em Matem\'{a}tica e Aplica\c{c}\~{o}%
es (CIMA),}\\
{\small \ Instituto de Investiga\c{c}\~{a}o e Forma\c{c}\~{a}o Avan\c{c}ada, 
}\\
{\small \ Universidade de \'{E}vora. Rua Rom\~{a}o Ramalho, 59, }\\
{\small \ 7000-671 \'{E}vora, Portugal}\\
$^{(\ddag )}${\small Faculdade de Ci\^{e}ncias e Tecnologia,}\\
{\small N$\acute{u}$cleo de Matem\'{a}tica e Aplica\c{c}\~{o}es (NUMAT),}\\
{\small Universidade de Cabo Verde. Campus de Palmarejo,}\\
{\small 279 Praia, Cabo Verde}}
\date{}
\maketitle

\begin{abstract}
This paper presents sufficient conditions for the solvability of the third
order three point boundary value problem%
\begin{equation*}
\left\{ 
\begin{array}{c}
-u^{\prime \prime \prime }(t)=f(t,\,v(t),\,v^{\prime }(t)) \\ 
-v^{\prime \prime \prime }(t)=h(t,\,u(t),\,u^{\prime }(t)) \\ 
u(0)=u^{\prime }(0)=0,u^{\prime }(1)=\alpha u^{\prime }(\eta ) \\ 
v(0)=v^{\prime }(0)=0,v^{\prime }(1)=\alpha v^{\prime }(\eta ).%
\end{array}%
\right.
\end{equation*}%
The arguments apply Green's function associated to the linear problem and
the Guo--Krasnosel'ski\u{\i} theorem of compression-expansion cones. The
dependence on the first derivatives is overcome by the construction of an
adequate cone and suitable conditions of superlinearity/sublinearity near $0$
and $+\infty .$ Last section contains an example to illustrate the
applicability of the theorem.\bigskip

\bigskip

\end{abstract}

\section{Introduction}

The solvability of systems of differential equations of second and higher
order, with different types of boundary conditions has received an
increasing interest in last years. See, for instance, \cite{Asif, Cui, Hend,
Hend2, Inf+FM+PP, Inf+PP, Kang, Lee, Li, Liu} and references therein .
However systems where the nonlinearities can depend on the first derivatives
are scarce (see \cite{Jank} ). Motivated by the works referred above, this
paper contributes to fill that gap, applying cones theory to the third order
three point boundary value problem%
\begin{equation}
\left\{ 
\begin{array}{c}
-u^{\prime \prime \prime }(t)=f(t,\,v(t),\,v^{\prime }(t)) \\ 
-v^{\prime \prime \prime }(t)=h(t,\,u(t),\,u^{\prime }(t)) \\ 
u(0)=u^{\prime }(0)=0,u^{\prime }(1)=\alpha u^{\prime }(\eta ) \\ 
v(0)=v^{\prime }(0)=0,v^{\prime }(1)=\alpha v^{\prime }(\eta ).%
\end{array}%
\right.  \label{eq1}
\end{equation}

The non-negative continuous functions $f,\;h\in C\left( [0,\,1]\times
\lbrack 0,\,+\infty )^{2},\;[0,\,+\infty )\right) $ verifying adequate
superlinear and sublinear conditions$,$ $0<\eta <1$ and the parameter $%
\alpha $ such that $1<\alpha <\frac{1}{\eta }.$

Third order differential equations can model various phenomena in physics,
biology or physiology such as the flow of a thin film of viscous fluid over
a solid surface (see\cite{Bernis, Tuck}), the solitary waves solution of the
Korteweg--de Vries equation (\cite{Liu+Chen}), or the thyroid-pituitary
interaction (\cite{Danzinger}).

A key point in our method is the fact that the Green's function associated
to the linear problem and its first derivative are nonnegative and verify
some adequate estimates. The existence of a positive and increasing solution
of the system (\ref{eq1}), is obtained by the well-known Guo--Krasnosel'ski%
\u{\i} theorem on cones compression-expansion. The dependence on the first
derivatives is overcome by the construction of an adequate cone and suitable
conditions of superlinearity/sublinearity near $0$ and $+\infty .$

The paper is organized in the following way: In section 2 we present the \ \
\ \ \ integral equations equivalent to problem (\ref{eq1}), the explicit
form of the Green's function, and its derivative, and the definition of some
functions used for its estimation. Section 3 contains the growth assumption
on the nonlinearities and the main result to prove the existence of an
increasing solution via cones theory. In last section an example illustrates
the applicability of the theorem.

\section{Preliminary results}

The pair of functions $(u(t),\;v(t))\in \big(C^{3}[0,\,1],\,(0,\;+\infty )%
\big)^{2}$ is a solution of problem (\ref{eq1}) if and only if $%
(u(t),\;v(t))\in \big(C^{3}[0,\,1],\,(0,\;+\infty )\big)^{2}$ it is a
solution of the following system of integral equations%
\begin{equation}
\left\{ 
\begin{array}{c}
u(t)=\int_{0}^{1}G(t,s)f(s,\,v(s),\,v^{\prime }(s))ds \\ 
\\ 
v(t)=\int_{0}^{1}G(t,s)h(s,\,u(s),\,u^{\prime }(s))ds,%
\end{array}%
\right.  \label{eq2}
\end{equation}%
where $G(t,s)$ is the Green's function associated to problem (\ref{eq1}),
defined by

\begin{equation}
G(t,s)=\frac{1}{2(1-\alpha \eta )}\left\{ 
\begin{array}{cc}
(2ts-s^{2})(1-\alpha \eta )+t^{2}s(\alpha -1) & s\leq \min \{\eta ,t\}, \\ 
t^{2}(1-\alpha \eta )+t^{2}s(\alpha -1) & t\leq s\leq \eta , \\ 
(2ts-s^{2})(1-\alpha \eta )+t^{2}(\alpha \eta -s) & \eta \leq s\leq t, \\ 
t^{2}(1-s) & \max \{\eta ,t\}\leq s.%
\end{array}%
\right.  \label{Green}
\end{equation}

Next Lemmas provide some properties of the Green's functions and its
derivative.

\begin{lemma}
(\cite{Li-jun})\label{Lema g0}Let $0<\eta <1$ and $1<\alpha <\frac{1}{\eta }$%
. Then for any $(t,s)\in \lbrack 0,\,1]\times \lbrack 0,\,1]$, we have $%
0\leq G(t,s)\leq g_{0}(s)$, where 
\begin{equation*}
g_{0}(s)=\frac{1+\alpha }{1-\alpha \eta }s(1-s).
\end{equation*}
\end{lemma}

\begin{lemma}
(\cite{Li-jun})\label{Lema k0}Let $0<\eta <1$ and $1<\alpha <\frac{1}{\eta }$%
. Then for any $(t,s)\in \lbrack \frac{\eta }{\alpha },\,\eta ]\times
\lbrack 0,\,1]$, the Green function $G(t,s)$ verifies $G(t,s)\geq
k_{0}g_{0}(s)$, where 
\begin{equation}
0<k_{0}:=\frac{\eta ^{2}}{2\alpha ^{2}(1+\alpha )}\min \{\alpha -1,\,1\}<1.
\label{k0}
\end{equation}
\end{lemma}

The derivative of $G$ is given by

\begin{equation*}
\frac{\partial G}{\partial t}(t,s)=\frac{1}{(1-\alpha \eta )}\left\{ 
\begin{array}{cc}
s(1-\alpha \eta )+ts(\alpha -1) & s\leq \min \{\eta ,t\}, \\ 
t(1-\alpha \eta )+ts(\alpha -1) & t\leq s\leq \eta , \\ 
s(1-\alpha \eta )+t(\alpha \eta -s) & \eta \leq s\leq t, \\ 
t(1-s) & \max \{\eta ,t\}\leq s,%
\end{array}%
\right.
\end{equation*}%
and verifies the following lemmas:

\begin{lemma}
\label{Lema g1}For $0<\eta <1,$ $1<\alpha <\frac{1}{\eta }$ and any $%
(t,s)\in \lbrack 0,\,1]\times \lbrack 0,\,1]$, we have $0\leq \frac{\partial
G}{\partial t}(t,s)\leq g_{1}(s)$, where 
\begin{equation*}
g_{1}(s)=\frac{(1-s)}{(1-\alpha \eta )}.
\end{equation*}
\end{lemma}

\begin{proof}
For $s\leq \min \{\eta ,t\}$, we have 
\begin{eqnarray*}
\frac{t(1-\alpha \eta )+ts(\alpha -1)}{(1-\alpha \eta )} &\leq &\frac{%
s(1-\alpha \eta )+s(\alpha -1)}{(1-\alpha \eta )}=\frac{s(\alpha -\alpha
\eta )}{(1-\alpha \eta )} \\
&=&\frac{s\alpha (1-\eta )}{(1-\alpha \eta )}\leq \frac{s\alpha (1-s)}{%
(1-\alpha \eta )}\leq \frac{(1-s)}{(1-\alpha \eta )}.
\end{eqnarray*}

If $t\leq s\leq \eta $, 
\begin{equation*}
\frac{t(1-\alpha \eta )+ts(\alpha -1)}{(1-\alpha \eta )}=\frac{t(1-\alpha
\eta +s\alpha -s)}{(1-\alpha \eta )}\leq \frac{(1-\alpha \eta +\eta \alpha
-s)}{(1-\alpha \eta )}=\frac{(1-s)}{(1-\alpha \eta )}.
\end{equation*}

For $\eta \leq s\leq t$, 
\begin{eqnarray*}
\frac{s(1-\alpha \eta )+t(\alpha \eta -s)}{(1-\alpha \eta )} &\leq &\frac{%
s(1-\alpha \eta )+(\alpha \eta -s)}{(1-\alpha \eta )}=\frac{\alpha \eta (1-s)%
}{(1-\alpha \eta )} \\
&\leq &\frac{\alpha s(1-s)}{(1-\alpha \eta )}\leq \frac{(1-s)}{(1-\alpha
\eta )}.
\end{eqnarray*}

If $\max \{\eta ,t\}\leq s$, then%
\begin{equation*}
\frac{t(1-s)}{(1-\alpha \eta )}\leq \frac{s(1-s)}{(1-\alpha \eta )}\leq 
\frac{(1-s)}{(1-\alpha \eta )}.
\end{equation*}
So, 
\begin{equation*}
\frac{\partial G}{\partial t}(t,s)\leq g_{1}(s):=\frac{(1-s)}{(1-\alpha \eta
)},\text{ for }(t,s)\in \lbrack 0,\,1]\times \lbrack 0,\,1].
\end{equation*}
\end{proof}

\begin{lemma}
\label{Lema k1}For $0<\eta <1,$ $1<\alpha <\frac{1}{\eta }$ and any $%
(t,s)\in \lbrack \frac{\eta }{\alpha },\,\eta ]\times \lbrack 0,\,1]$, the
derivative of the Green function $\frac{\partial G}{\partial t}(t,s)$
verifies $\frac{\partial G}{\partial t}(t,s)\geq k_{1}g_{1}(s)$, with 
\begin{equation}
0<k_{1}:=\min \{\alpha \eta ,\;\eta \}<1.  \label{k1}
\end{equation}
\end{lemma}

\begin{proof}
To find $k_{1}$ such that%
\begin{equation*}
k_{1}g_{1}(s)\leq \frac{\partial G}{\partial t}(t,s),
\end{equation*}%
we evaluate it in each branch of $\frac{\partial G}{\partial t}(t,s)$ for $%
(t,s)\in \lbrack \frac{\eta }{\alpha },\,\eta ]\times \lbrack 0,\,1]$,

\begin{description}
\item[(i)] For $s\leq \min \{\eta ,t\}$, we must have%
\begin{equation*}
k_{1}\frac{1-s}{1-\alpha \eta }\leq \frac{s(1-\alpha \eta )+ts(\alpha -1)}{%
1-\alpha \eta },
\end{equation*}%
that is%
\begin{equation*}
k_{1}\leq \frac{s(1-\alpha \eta )+ts(\alpha -1)}{1-s}\leq \frac{\eta
(1-\alpha \eta )+\eta (\alpha -1)}{1-\eta }\leq \frac{\eta (\alpha -\alpha
\eta )}{1-\eta }=\alpha \eta <1.
\end{equation*}

\item[(ii)] If $t\leq s\leq \eta $, the inequality 
\begin{equation*}
k_{1}\frac{1-s}{1-\alpha \eta }\leq \frac{t(1-\alpha \eta +s\alpha -s)}{%
1-\alpha \eta }
\end{equation*}%
holds for 
\begin{equation*}
k_{1}\leq \frac{t(1-\alpha \eta +s\alpha -s)}{1-s}\leq \frac{\eta (1-\alpha
\eta +\eta \alpha -\eta )}{1-\eta }\leq \frac{\eta (1-\eta )}{1-\eta }=\eta
<1.
\end{equation*}%
So, we can to define%
\begin{equation*}
0<k_{1}=\min \{\alpha \eta ,\;\eta \}<1.
\end{equation*}
\end{description}
\end{proof}

The existence tool will be the well known Guo-Krasnoselskii results in
expansive and compressive cones theory:

\begin{lemma}
(\cite{Guo})\label{Lema G-K}Let $(E,\,\Vert \cdot \Vert )$ be a Banach
space, and $P\subset E$ be a cone in $E$. Assume that $\Omega _{1}$ and $%
\Omega _{2}$ are open subsets of $E$ such that $0\in \Omega _{1},\;\overline{%
\Omega _{1}}\subset \Omega _{2}$.

If $T\,:\,P\cap (\overline{\Omega _{2}}\setminus \Omega _{1})\,\rightarrow
\,P$ is a completely continuous operator such that either

\begin{itemize}
\item[(i)] $\Vert Tu\Vert \leq \Vert u\Vert ,\;u\in P\cap \partial \Omega
_{1},$ and $\Vert Tu\Vert \geq \Vert u\Vert ,\;u\in P\cap \partial \Omega
_{2}$,

or

\item[(ii)] $\Vert Tu\Vert \geq \Vert u\Vert ,\;u\in P\cap \partial \Omega
_{1},$ and $\Vert Tu\Vert \leq \Vert u\Vert ,\;u\in P\cap \partial \Omega
_{2}$,

\noindent then $T$ has a fixed point in $P\cap (\overline{\Omega _{2}}%
\backslash \Omega _{1})$.
\end{itemize}
\end{lemma}

\section{Main result}

Consider the following growth assumptions 
\begin{align*}
(A1)\;\;\ \ \ \underset{t\in \lbrack 0,1],\text{ }\Vert v\Vert
_{C^{1}}\rightarrow 0\;}{\limsup }\frac{f(t,\,v,\,v^{\prime })}{%
|v|+|v^{\prime }|}& =0\text{ \ and \ }\underset{t\in \lbrack 0,1],\text{ }%
\Vert u\Vert _{C^{1}}\rightarrow 0\;}{\limsup }\frac{h(t,\,u,\,u^{\prime })}{%
|u|+|u^{\prime }|}=0; \\
(A2)\;\;\underset{t\in \lbrack 0,1],\text{ }\Vert v\Vert
_{_{C^{1}}}\rightarrow +\infty }{\liminf }\frac{f(t,\,v,\,v^{\prime })}{%
|v|+|v^{\prime }|}& =+\infty \text{ \ and \ }\underset{t\in \lbrack 0,1],%
\text{ }\Vert u\Vert _{_{C^{1}}}\rightarrow +\infty }{\liminf }\frac{%
h(t,\,u,\,u^{\prime })}{|u|+|u^{\prime }|}=+\infty ; \\
(A3)\;\ \;\ \ \underset{t\in \lbrack 0,1],\text{ }\Vert v\Vert
_{C^{1}}\rightarrow 0}{\liminf }\frac{f(t,\,v,\,v^{\prime })}{|v|+|v^{\prime
}|}& =+\infty \text{ \ and \ }\underset{t\in \lbrack 0,1],\text{ }\Vert
u\Vert _{C^{1}}\rightarrow 0}{\liminf }\frac{h(t,\,u,\,u^{\prime })}{%
|u|+|u^{\prime }|}=+\infty ; \\
(A4)\;\;\underset{t\in \lbrack 0,1],\text{ }\Vert v\Vert
_{_{C^{1}}}\rightarrow +\infty }{\limsup }\frac{f(t,\,v,\,v^{\prime })}{%
|v|+|v^{\prime }|}& =0\text{ \ and \ }\underset{t\in \lbrack 0,1],\text{ }%
\Vert u\Vert _{_{C^{1}}}\rightarrow +\infty }{\limsup }\frac{%
h(t,\,u,\,u^{\prime })}{|u|+|u^{\prime }|}=0.
\end{align*}

The main result is given by next theorem :

\begin{theorem}
\label{Main thm}Let $f,\;h:[0,\,1]\times \lbrack 0,\,+\infty
)^{2}\rightarrow \lbrack 0,\,+\infty )$ be continuous functions such that
assumptions $(A1)$ and $(A2),$ or $(A3)$ and $(A4),$ hold.\newline
Then problem (\ref{eq1}) has at least one positive solution $%
(u(t),\,v(t))\in \big(C^{3}[0,\,1]\big)^{2}$, that is $u(t)>0$, $v(t)>0,$ $%
\forall t\in \lbrack 0,\,1]$.
\end{theorem}

\begin{proof}
Let $E=C^{1}[0,\,1]$ be the Banach space equipped with the norm $\Vert \cdot
\Vert _{C^{1}}$, defined by $\Vert w\Vert _{C^{1}}:=\max \left\{ \Vert
w\Vert ,\Vert w^{\prime }\Vert \right\} $ and $\Vert y\Vert :=\underset{t\in
\lbrack 0,\,1]\,}{\max }|y(t)|$.

Consider the set 
\begin{equation*}
K=\left\{ w\in E\,:\,w(t)\geq 0,\,\underset{t\in \left[ \frac{\eta }{\alpha }%
,\,\eta \right] \,}{\min }w(t)\geq k_{0}\Vert w\Vert ,\text{ }\underset{t\in %
\left[ \frac{\eta }{\alpha },\,\eta \right] \,}{\min }w^{\prime }(t)\geq
k_{1}\Vert w^{\prime }\Vert \right\} ,
\end{equation*}%
with $k_{0}$ and $k_{1}$ given by (\ref{k0}) and (\ref{k1}), respectively,
and the operators $T_{1}\,:\,K\,\rightarrow K$ $\ $and $T_{2}\,:\,K\,%
\rightarrow \,K$ such that%
\begin{equation}
\left\{ 
\begin{array}{c}
T_{1}u(t)=\int_{0}^{1}G(t,s)f(s,\,v(s),\,v^{\prime }(s))ds \\ 
\\ 
T_{2}v(t)=\int_{0}^{1}G(t,s)h(s,\,u(s),\,u^{\prime }(s))ds.%
\end{array}%
\right.  \label{DefT1T2}
\end{equation}

By (\ref{eq2}), the solutions of the initial system (\ref{eq1}) are fixed
points of the operator $T:=(T_{1},T_{2}).$

First we show that $K$ is a cone. By definition of $K$ it is clear that $K$
is not identically zero or empty.

Consider $a,\;b\in \mathbb{R}^{+}$ and $\forall x,\;y\in K$. Then%
\begin{equation*}
x\in K\Rightarrow x\in E\,:\,x(t)\geq 0,\,\underset{t\in \lbrack 0,\,1]\,}{%
\min }x(t)\geq k_{0}\Vert x\Vert ,\text{ }\underset{t\in \lbrack 0,\,1]\,}{%
\min }x^{\prime }(t)\geq k_{1}\Vert x^{\prime }\Vert ,
\end{equation*}%
\begin{equation*}
y\in K\Rightarrow y\in E\,:\,y(t)\geq 0,\,\underset{t\in \lbrack 0,\,1]\,}{%
\min }y(t)\geq k\Vert y\Vert ,\,\underset{t\in \lbrack 0,\,1]\,}{\min }%
y^{\prime }(t)\geq k_{1}\Vert y^{\prime }\Vert .
\end{equation*}

As $E$ is a vector space, consider the linear combination $ax+by\in E.$ 
\begin{eqnarray*}
\underset{t\in \lbrack 0,\,1]\,}{\min }\left( ax(t)+by(t)\right) &=&a%
\underset{t\in \lbrack 0,\,1]\,}{\min }x(t)+b\underset{t\in \lbrack 0,\,1]\,}%
{\min }y(t) \\
&\geq &ak_{0}\Vert x\Vert +bk_{0}\Vert y\Vert =k_{0}\left( a\Vert x\Vert
+b\Vert y\Vert \right) \\
&\geq &k_{0}\left\Vert ax(t)+by(t)\right\Vert ,
\end{eqnarray*}%
and%
\begin{eqnarray*}
\underset{t\in \lbrack 0,\,1]\,}{\min }\left( ax(t)+by(t)\right) ^{\prime }
&=&a\underset{t\in \lbrack 0,\,1]\,}{\min }\left( x(t)\right) ^{\prime }+b%
\underset{t\in \lbrack 0,\,1]\,}{\min }\left( y(t)\right) ^{\prime } \\
&\geq &ak_{1}\Vert x^{\prime }\Vert +bk_{1}\Vert y^{\prime }\Vert
=k_{1}\left( a\Vert x^{\prime }\Vert +b\Vert y^{\prime }\Vert \right) \\
&\geq &k_{0}\left\Vert \left( ax(t)+by(t)\right) ^{\prime }\right\Vert .
\end{eqnarray*}%
Therefore $ax+by\in K,$ that is $K$ is a cone.

Now we show that $T_{1}$ and $T_{2}$ are completely continuous, i.e, are
equicontinuous and uniformly bounded.

For the reader's convenience the proof for $T_{1}$ will follow several steps
and claims. The arguments for $T_{2}$ are analogous.\medskip

\textbf{Step 1:} $T_{1}$\textit{\ and }$T_{2}$\textit{\ are well defined in }%
$K.$\medskip

To prove that $T_{1}K\subset K$ consider $u\in K.$

As $G(t,s)\geq 0$ for $(t,s)\in \lbrack 0,\,1]\times \lbrack 0,\,1],$ It is
clear that $T_{1}u(t)\geq 0.$

By Lemma \ref{Lema g0}, the positivity of $f$ and (\ref{DefT1T2}),%
\begin{equation*}
0\leq T_{1}u(t)=\int_{0}^{1}G(t,s)f(s,\,v(s),\,v^{\prime }(s))ds\leq
\int_{0}^{1}g_{0}(s)f(s,\,v(s),\,v^{\prime }(s))ds.
\end{equation*}

So,%
\begin{equation}
\Vert T_{1}u\Vert \leq \int_{0}^{1}g_{0}(s)f(s,\,v(s),\,v^{\prime }(s))ds.
\label{eq4}
\end{equation}%
From Lemma \ref{Lema k0} and (\ref{eq4}),%
\begin{eqnarray*}
T_{1}u(t) &=&\int_{0}^{1}G(t,s)f(s,\,v(s),\,v^{\prime }(s))ds \\
&\geq &k_{0}\int_{0}^{1}g_{0}(s)f(s,\,v(s),\,v^{\prime }(s))ds\geq
k_{0}\Vert T_{1}u\Vert ,\;\text{for }t\in \left[ \frac{\eta }{\alpha },\eta %
\right] ,
\end{eqnarray*}%
with $k_{0}$ given by (\ref{k0}). By Lemma \ref{Lema g1},%
\begin{equation*}
(T_{1}u(t))^{\prime }=\int_{0}^{1}\frac{\partial G}{\partial t}%
(t,s)f(s,\,v(s),\,v^{\prime }(s))ds\leq
\int_{0}^{1}g_{1}(s)f(s,\,v(s),\,v^{\prime }(s))ds,
\end{equation*}%
So,%
\begin{equation}
\Vert (T_{1}u)^{\prime }\Vert \leq
\int_{0}^{1}g_{1}(s)f(s,\,v(s),\,v^{\prime }(s))ds.  \label{eq5}
\end{equation}

By Lemma \ref{Lema k1} and (\ref{eq5}), it follows 
\begin{eqnarray*}
(T_{1}u(t))^{\prime } &=&\int_{0}^{1}\frac{\partial G}{\partial t}%
(t,s)f(s,\,v(s),\,v^{\prime }(s))ds\geq
k_{1}\int_{0}^{1}g_{1}(s)f(s,\,v(s),\,v^{\prime }(s)ds \\
&\geq &k_{1}\Vert (T_{1}u)^{\prime }\Vert ,\;\text{for }t\in \left[ \frac{%
\eta }{\alpha },\eta \right] ,
\end{eqnarray*}%
and $k_{1}$ as in (\ref{k1}).

So $T_{1}K\subset K.$ Analogously it can be shown that $T_{2}K\subset
K.\medskip $

Assume that $(A1)$ and $(A2)$ hold.

By $(A1),$ there exists $0<\delta _{1}<1$ such that, for$(t,v,v^{\prime
})\in \lbrack 0,1]\times \lbrack 0,\delta _{1}]^{2}$ and $(t,u,u^{\prime
})\in \lbrack 0,1]\times \lbrack 0,\delta _{1}]^{2},$ 
\begin{equation}
f(s,\,v(s),\,v^{\prime }(s))\leq \varepsilon _{1}\left( |v(s)|+|v^{\prime
}(s)|\right)  \label{e1}
\end{equation}%
and 
\begin{equation}
h(s,\,u(s),\,u^{\prime }(s))\leq \varepsilon _{2}\left( |u(s)|+|u^{\prime
}(s)|\right) ,  \label{e2}
\end{equation}%
with $\varepsilon _{1}$ and $\varepsilon _{2}$ to be defined forward.\medskip

\textbf{Step 2: \ }$T_{1}$\textit{and }$T_{2}$\textit{\ are completely
continuous in }$C^{1}[0,\,1]$\textit{.}\medskip

$T_{1}$ is continuous in $C^{1}[0,\,1]$ as: $G(t,s),$ $\frac{\partial G}{%
\partial t}(t,s)$ and $f$ are continuous.

Consider the set $B\subset K,$ bounded in $C^{1},$ and let $u,v\in B$. Then
there are $M_{1},M_{2}>0$ such that $\Vert u\Vert _{C^{1}}<M_{1}$ and $\Vert
v\Vert _{C^{1}}<M_{2}.$\medskip

\textbf{Claim 2.1}. $T_{1}$\textit{\ is uniformly bounded in }$C^{1}[0,\,1]$%
. \medskip

In fact, by (\ref{e1}), there are $M_{3},M_{4}>0$ such that%
\begin{eqnarray*}
\Vert {T_{1}}u\Vert &=&\underset{t\in \lbrack 0,\,1]}{\max }\left\vert {T_{1}%
}u(t)\right\vert =\underset{t\in \lbrack 0,\,1]}{\max }\left\vert
\int_{0}^{1}G(t,s)f(s,\,v(s),\,v^{\prime }(s))ds\right\vert \\
&\leq &\int_{0}^{1}\underset{t\in \lbrack 0,\,1]}{\max }\left\vert
G(t,s)\right\vert \left\vert f(s,\,v(s),\,v^{\prime }(s))\right\vert ds \\
&\leq &\int_{0}^{1}\underset{t\in \lbrack 0,\,1]}{\max }\left\vert
G(t,s)\right\vert \varepsilon _{1}\left( |v(s)|+|v^{\prime }(s)|\right) ds \\
&\leq &2\varepsilon _{1}\Vert v\Vert _{C^{1}}\int_{0}^{1}\underset{t\in
\lbrack 0,\,1]}{\max }\left\vert G(t,s)\right\vert <M_{3},\;\forall u\in B,
\end{eqnarray*}%
\begin{eqnarray*}
\Vert \left( {T_{1}}u\right) ^{\prime }\Vert &=&\underset{t\in \lbrack 0,\,1]%
}{\max }\left\vert \left( {T_{1}}u(t)\right) ^{\prime }\right\vert =\underset%
{t\in \lbrack 0,\,1]}{\max }\left\vert \int_{0}^{1}\frac{\partial G}{%
\partial t}(t,s)f(s,\,v(s),\,v^{\prime }(s)ds\right\vert \\
&\leq &\int_{0}^{1}\underset{t\in \lbrack 0,\,1]}{\max }\left\vert \frac{%
\partial G}{\partial t}(t,s)\right\vert \left\vert f(s,\,v(s),\,v^{\prime
}(s)\right\vert ds \\
&\leq &\int_{0}^{1}\underset{t\in \lbrack 0,\,1]}{\max }\left\vert \frac{%
\partial G}{\partial t}(t,s)\right\vert \varepsilon _{1}\left(
|v(s)|+|v^{\prime }(s)|\right) ds \\
&\leq &2\varepsilon _{1}\Vert v\Vert _{C^{1}}\int_{0}^{1}\underset{t\in
\lbrack 0,\,1]}{\max }\left\vert \frac{\partial G}{\partial t}%
(t,s)\right\vert ds<M_{4},\;\forall u\in B.
\end{eqnarray*}

Defining $M:=\max \left\{ M_{3},M_{4}\right\} ,$ \ then $\Vert {T_{1}}u\Vert
_{C^{1}}\leq M$.\medskip

\textbf{Claim 2.2}. $T_{1}$\textit{\ is equicontinuous in }$C^{1}[0,\,1]$.
\medskip

Let $t_{1}$ and $t_{2}\in \lbrack 0,1].$ Without loss of generality suppose $%
t_{1}\leq t_{2}.$ So%
\begin{eqnarray*}
\left\vert Tu(t_{1})-Tu(t_{2})\right\vert &=&\left\vert \int_{0}^{1}\left[
G(t_{1},s)-G(t_{2},s)\right] f(s,\,v(s),\,v^{\prime }(s))ds\right\vert \\
&\leq &\int_{0}^{1}\left\vert G(t_{1},s)-G(t_{2},s)\right\vert \varepsilon
_{1}\left( |v|+|v^{\prime }|\right) ds \\
&\leq &2\varepsilon _{1}\Vert v\Vert _{C^{1}}\int_{0}^{1}\left\vert
G(t_{1},s)-G(t_{2},s)\right\vert ds\rightarrow 0,\;\text{as}t_{1}\rightarrow
t_{2},\;
\end{eqnarray*}%
and%
\begin{eqnarray*}
\left\vert \left( Tu(t_{1})\right) ^{\prime }-\left( Tu(t_{2})\right)
^{\prime }\right\vert &=&\left\vert \int_{0}^{1}\left[ \frac{\partial G}{%
\partial t}(t_{1},s)-\frac{\partial G}{\partial t}(t_{2},s)\right]
f(s,\,v(s),\,v^{\prime }(s))ds\right\vert \\
&\leq &\int_{0}^{1}\left\vert \frac{\partial G}{\partial t}(t_{1},s)-\frac{%
\partial G}{\partial t}(t_{2},s)\right\vert \varepsilon _{1}\left(
|v|+|v^{\prime }|\right) ds \\
&\leq &2\varepsilon _{1}\Vert v\Vert _{C^{1}}\int_{0}^{1}\left\vert \frac{%
\partial G}{\partial t}(t_{1},s)-\frac{\partial G}{\partial t}%
(t_{2},s)\right\vert ds\rightarrow 0,\text{ as }t_{1}\rightarrow t_{2}.
\end{eqnarray*}

By the Arz\`{e}la-Ascoli's lemma $T_{1}B$ is relatively compact, that is, $%
T_{1}$ is compact .

Applying the same technique, using (\ref{e2}), it can be shown that $T_{2}$
is compact, too. Consequently $T$ is compact. \medskip

Next steps will prove that assumptions of Lemma \ref{Lema G-K} hold.\medskip

\textbf{Step 3: \ }$\Vert T_{1}u\Vert _{C^{1}}\leq \Vert u\Vert _{C^{1}}$%
\textit{, for some }$\rho _{1}>0$ \textit{and} \textit{\ }$u\in K\cap
\partial \Omega _{1}$\textit{\ with }$\Omega _{1}=\{u\in E\,:\,\Vert u\Vert
_{C^{1}}<\rho _{1}\}.$\medskip

By $(A1),$ define $0<\rho _{1}<1$ such that $(t,\,v,v^{\prime })\in \lbrack
0,\,1]\times \lbrack 0,\,\rho _{1}]^{2}$ and $(t,\,u,u^{\prime })\in \lbrack
0,\,1]\times \lbrack 0,\,\rho _{1}]^{2}$.

From (\ref{e1}) and (\ref{e2}), choose $\varepsilon _{1},\varepsilon _{2}>0$
sufficiently small such that%
\begin{equation}
\max \text{ }\left\{ 
\begin{array}{c}
\varepsilon _{1}\varepsilon _{2}\int_{0}^{1}g_{0}(s)ds\int_{0}^{1}\left(
g_{0}(r)+g_{1}(r)\right) dr,\text{ } \\ 
\varepsilon _{1}\varepsilon _{2}\int_{0}^{1}g_{1}(s)ds\int_{0}^{1}\left(
g_{0}(r)+g_{1}(r)\right) dr%
\end{array}%
\right\} <\frac{1}{2}.  \label{eps}
\end{equation}

If $u\in K$ and $\Vert u\Vert _{C^{1}}=\rho _{1}$, then, by Lemma \ref{Lema
g0}, (\ref{eq2}) and (\ref{eps}), 
\begin{eqnarray*}
T_{1}u(t) &\leq &\int_{0}^{1}g_{0}(s)\varepsilon _{1}\left( |v|+|v^{\prime
}|\right) ds \\
&\leq &\int_{0}^{1}g_{0}(s)\varepsilon _{1}\int_{0}^{1}\left( \left\vert
G(t,r)\right\vert +\left\vert \frac{\partial G}{\partial t}(t,r)\right\vert
\right) \left\vert h(r,\,u(r),\,u^{\prime }(r))\right\vert drds \\
&\leq &\int_{0}^{1}g_{0}(s)\varepsilon _{1}\int_{0}^{1}\left(
g_{0}(r)+g_{1}(r)\right) \left\vert h(r,\,u(r),\,u^{\prime }(r))\right\vert
drds \\
&\leq &\varepsilon _{1}\varepsilon
_{2}\int_{0}^{1}g_{0}(s)ds\int_{0}^{1}\left( g_{0}(r)+g_{1}(r)\right) \left(
|u(r)|+|u^{\prime }(r)|\right) dr \\
&\leq &2\varepsilon _{1}\varepsilon _{2}\Vert u\Vert
_{C^{1}}\int_{0}^{1}g_{0}(s)ds\int_{0}^{1}\left( g_{0}(r)+g_{1}(r)\right)
dr<\Vert u\Vert _{C^{1}},
\end{eqnarray*}%
and%
\begin{eqnarray*}
\left( T_{1}u(t)\right) ^{\prime } &=&\int_{0}^{1}\frac{\partial G}{\partial
t}(t,s)f(s,\,v(s),\,v^{\prime }(s)ds\leq \int_{0}^{1}g_{1}(s)\varepsilon
_{1}\left( |v|+|v^{\prime }|\right) ds \\
&\leq &\int_{0}^{1}g_{1}(s)\varepsilon _{1}\int_{0}^{1}\left(
g_{0}(r)+g_{1}(r)\right) \left\vert h(r,\,u(r),\,u^{\prime }(r))\right\vert
drds \\
&\leq &2\varepsilon _{1}\varepsilon _{2}\Vert u\Vert
_{C^{1}}\int_{0}^{1}g_{1}(s)ds\int_{0}^{1}\left( g_{0}(r)+g_{1}(r)\right)
dr<\Vert u\Vert _{C^{1}},
\end{eqnarray*}

Therefore $\Vert T_{1}u\Vert _{C^{1}}\leq \Vert u\Vert _{C^{1}}.$\medskip

\textbf{Step 4: \ }$\Vert T_{1}u\Vert _{C^{1}}\geq \Vert u\Vert _{C^{1}}$%
\textit{, for some }$\rho _{2}>0$ \textit{and} \textit{\ }$u\in K\cap
\partial \Omega _{2}$\textit{\ with }$\Omega _{2}=\{u\in E\,:\,\Vert u\Vert
_{C^{1}}<\rho _{2}\}.$\medskip

By $(A2),$ \ $\Vert v\Vert _{C^{1}}\rightarrow +\infty $ \ and $\Vert u\Vert
_{C^{1}}\rightarrow +\infty .$ Therefore there are several cases to be
considered:\medskip

\textbf{Case 4.1.} \textit{Suppose that there exist }$\theta _{1},\theta
_{2}>0$\textit{\ such that }$\Vert v\Vert \rightarrow +\infty ,$\textit{\ }$%
\Vert v^{\prime }\Vert \leq $\textit{\ }$\theta _{1},$\textit{\ }$\Vert
u\Vert \rightarrow +\infty $\textit{\ and }$\Vert u^{\prime }\Vert \leq $%
\textit{\ }$\theta _{2}.$\medskip

Consider $\rho >0$ such that for $(t,v,v^{\prime })\in \lbrack 0,\,1]\times
\lbrack \rho ,\,+\infty )\times \left[ 0,\theta _{1}\right] $ and $%
(t,\,u,u^{\prime })\in \lbrack 0,\,1]\times \lbrack \rho ,\,+\infty )\times %
\left[ 0,\theta _{2}\right] $, we have%
\begin{equation}
f(t,v(t),v^{\prime }(t))\geq \xi _{1}\left( |v(t)|+|v^{\prime }(t)|\right)
\label{csi1}
\end{equation}%
and%
\begin{equation}
h(t,u(t),u^{\prime }(t))\geq \xi _{2}\left( |u(t)|+|u^{\prime }(t)|\right) ,
\label{csi2}
\end{equation}%
with $\xi _{1}$, $\xi _{2}$ such that 
\begin{equation}
\min \left\{ 
\begin{array}{c}
\left( k_{0}\right) ^{2}\xi _{1}\xi _{2}\int_{\frac{\eta }{\alpha }}^{\eta
}g_{0}(s)ds\int_{\frac{\eta }{\alpha }}^{\eta }\left(
k_{0}g_{0}(r)+k_{1}g_{1}(r)\right) dr,\text{ } \\ 
\xi _{1}\xi _{2}k_{0}k_{1}\int_{\frac{\eta }{\alpha }}^{\eta
}g_{0}(s)ds\int_{\frac{\eta }{\alpha }}^{\eta }\left(
k_{0}g_{0}(r)+k_{1}g_{1}(r)\right) dr, \\ 
\xi _{1}\xi _{2}k_{0}k_{1}\int_{\frac{\eta }{\alpha }}^{\eta
}g_{1}(s)ds\int_{\frac{\eta }{\alpha }}^{\eta }\left(
k_{0}g_{0}(r)+k_{1}g_{1}(r)\right) dr,\text{ } \\ 
\left( k_{1}\right) ^{2}\xi _{1}\xi _{2}\int_{\frac{\eta }{\alpha }}^{\eta
}g_{1}(s)\text{ }ds\int_{\frac{\eta }{\alpha }}^{\eta }\left(
k_{0}g_{0}(r)+k_{1}g_{1}(r)\right) \text{ }dr \\ 
k_{0}\left( k_{0}+k_{1}\right) \xi _{1}\xi _{2}\int_{\frac{\eta }{\alpha }%
}^{\eta }g_{0}(s)\text{ }ds\int_{\frac{\eta }{\alpha }}^{\eta }\left(
k_{0}g_{0}(r)+k_{1}g_{1}(r)\right) \text{ }dr,\text{ } \\ 
k_{0}\left( k_{0}+k_{1}\right) \xi _{1}\xi _{2}\int_{\frac{\eta }{\alpha }%
}^{\eta }g_{1}(s)\text{ }ds\int_{\frac{\eta }{\alpha }}^{\eta }\left(
k_{0}g_{0}(r)+k_{1}g_{1}(r)\right) \text{ }dr%
\end{array}%
\right\} >1,  \label{eq7}
\end{equation}%
with $k_{0},k_{1}$ as in (\ref{k0}) and (\ref{k1}).

Let $u,v\in K$ such that $\Vert u\Vert _{C^{1}}=\rho _{2}$, where $\rho
_{2}:=\max \left\{ 2\rho _{1},\,\frac{\rho }{k_{0}},\,\frac{\rho }{k_{1}}%
\right\} $.

Then $\Vert u\Vert _{C^{1}}=\Vert u\Vert =\rho _{2}$ and $u(t)\geq
k_{0}\Vert u\Vert _{C^{1}}=k_{0}\rho _{2}\geq \rho ,\,t\in \lbrack 0,\,1].$
Similarly, $\Vert v\Vert _{C^{1}}=\Vert v\Vert =\rho _{2}$ and $v(t)\geq
k_{1}\Vert v\Vert _{C^{1}}=k_{0}\rho _{2}\geq \rho .$

By Lemma \ref{Lema k0}, (\ref{eq2}) and (\ref{eq7}),%
\begin{eqnarray*}
T_{1}u(t) &\geq &\int_{\frac{\eta }{\alpha }}^{\eta
}G(t,s)f(s,\,v(s),\,v^{\prime }(s))ds \\
&\geq &k_{0}\int_{\frac{\eta }{\alpha }}^{\eta }g_{0}(s)\text{ }%
f(s,\,v(s),\,v^{\prime }(s))ds\geq k_{0}\xi _{1}\int_{\frac{\eta }{\alpha }%
}^{\eta }g_{0}(s)\text{ }\left( |v(s)|+|v^{\prime }(s)|\right) ds \\
&=&k_{0}\xi _{1}\int_{\frac{\eta }{\alpha }}^{\eta }g_{0}(s)\text{ }ds\int_{%
\frac{\eta }{\alpha }}^{\eta }\left( \left\vert G(t,r)\right\vert
+\left\vert \frac{\partial G}{\partial t}(t,r)\right\vert \right) \left\vert
h(r,\,u(r),\,u^{\prime }(r))\right\vert dr \\
&\geq &k_{0}\xi _{1}\xi _{2}\int_{\frac{\eta }{\alpha }}^{\eta }g_{0}(s)%
\text{ }ds\int_{\frac{\eta }{\alpha }}^{\eta }\left(
k_{0}g_{0}(r)+k_{1}g_{1}(r)\right) \text{ }\left( |u(r)|+|u^{\prime
}(r)|\right) dr \\
&=&k_{0}\xi _{1}\int_{\frac{\eta }{\alpha }}^{\eta }g_{0}(s)\text{ }ds\int_{%
\frac{\eta }{\alpha }}^{\eta }\left( \left\vert G(t,r)\right\vert
+\left\vert \frac{\partial G}{\partial t}(t,r)\right\vert \right) \left\vert
h(r,\,u(r),\,u^{\prime }(r))\right\vert dr \\
&\geq &k_{0}\xi _{1}\xi _{2}\int_{\frac{\eta }{\alpha }}^{\eta }g_{0}(s)%
\text{ }ds\int_{\frac{\eta }{\alpha }}^{\eta }\left(
k_{0}g_{0}(r)+k_{1}g_{1}(r)\right) \text{ }\left( |u(r)|+|u^{\prime
}(r)|\right) dr \\
&\geq &k_{0}\xi _{1}\xi _{2}\int_{\frac{\eta }{\alpha }}^{\eta }g_{0}(s)%
\text{ }ds\int_{\frac{\eta }{\alpha }}^{\eta }\left(
k_{0}g_{0}(r)+k_{1}g_{1}(r)\right) \text{ }\left( \underset{r\in \left[ 
\frac{\eta }{\alpha },\,\eta \right] \,}{\min }u(r)+\underset{r\in \left[ 
\frac{\eta }{\alpha },\,\eta \right] \,}{\min }u^{\prime }(r)\right) dr \\
&\geq &k_{0}\xi _{1}\xi _{2}\int_{\frac{\eta }{\alpha }}^{\eta }g_{0}(s)%
\text{ }ds\int_{\frac{\eta }{\alpha }}^{\eta }\left(
k_{0}g_{0}(r)+k_{1}g_{1}(r)\right) \text{ }\left( k_{0}\Vert u\Vert
+k_{1}\Vert u^{\prime }\Vert \right) dr
\end{eqnarray*}%
\begin{eqnarray*}
&\geq &k_{0}\xi _{1}\xi _{2}\int_{\frac{\eta }{\alpha }}^{\eta }g_{0}(s)%
\text{ }ds\int_{\frac{\eta }{\alpha }}^{\eta }\left(
k_{0}g_{0}(r)+k_{1}g_{1}(r)\right) \text{ }k_{0}\Vert u\Vert _{C^{1}}dr \\
&=&k_{0}^{2}\Vert u\Vert _{C^{1}}\xi _{1}\xi _{2}\int_{\frac{\eta }{\alpha }%
}^{\eta }g_{0}(s)\text{ }ds\int_{\frac{\eta }{\alpha }}^{\eta }\left(
k_{0}g_{0}(r)+k_{1}g_{1}(r)\right) \text{ }dr>\Vert u\Vert _{C^{1}},
\end{eqnarray*}%
and, analogously,%
\begin{eqnarray*}
\left( T_{1}u(t)\right) ^{\prime } &\geq &\int_{\frac{\eta }{\alpha }}^{\eta
}\frac{\partial G}{\partial t}(t,s)f(s,\,v(s),\,v^{\prime }(s))ds \\
&\geq &k_{1}\int_{\frac{\eta }{\alpha }}^{\eta }g_{1}(s)\text{ }%
f(s,\,v(s),\,v^{\prime }(s))ds\geq k_{1}\xi _{1}\int_{\frac{\eta }{\alpha }%
}^{\eta }g_{1}(s)\text{ }\left( |v(s)|+|v^{\prime }(s)|\right) ds \\
&\geq &k_{1}\xi _{1}\xi _{2}\int_{\frac{\eta }{\alpha }}^{\eta }g_{1}(s)%
\text{ }ds\int_{\frac{\eta }{\alpha }}^{\eta }\left(
k_{0}g_{0}(r)+k_{1}g_{1}(r)\right) \text{ }k_{0}\Vert u\Vert _{C^{1}}dr \\
&=&k_{1}k_{0}\Vert u\Vert _{C^{1}}\xi _{1}\xi _{2}\int_{\frac{\eta }{\alpha }%
}^{\eta }g_{1}(s)\text{ }ds\int_{\frac{\eta }{\alpha }}^{\eta }\left(
k_{0}g_{0}(r)+k_{1}g_{1}(r)\right) \text{ }dr>\Vert u\Vert _{C^{1}}.
\end{eqnarray*}

Therefore $\Vert T_{1}u\Vert _{C^{1}}\geq \Vert u\Vert _{C^{1}}.$\medskip

\textbf{Case 4.2.} \textit{Suppose that there exist }$\theta _{3},\theta
_{4}>0$\textit{\ such that }$\Vert v^{\prime }\Vert \rightarrow +\infty ,$%
\textit{\ }$\Vert v\Vert \leq $\textit{\ }$\theta _{3},$\textit{\ }$\Vert
u^{\prime }\Vert \rightarrow +\infty $\textit{\ and }$\Vert u\Vert \leq $%
\textit{\ }$\theta _{4}.$\medskip

Consider $\rho >0$ such that for $(t,v,v^{\prime })\in \lbrack 0,\,1]\times %
\left[ 0,\theta _{3}\right] \times \lbrack \rho ,\,+\infty )$ and $%
(t,\,u,u^{\prime })\in \lbrack 0,\,1]\times \left[ 0,\theta _{4}\right]
\times \lbrack \rho ,\,+\infty )$, conditions (\ref{csi1}), (\ref{csi2}) and
(\ref{eq7}) hold.

Let $u,v\in K$ such that $\Vert u\Vert _{C^{1}}=\rho _{2}$, where $\rho
_{2}:=\max \left\{ 2\rho _{1},\,\frac{\rho }{k_{0}},\,\frac{\rho }{k_{1}}%
\right\} $.

Then $\Vert u\Vert _{C^{1}}=\Vert u^{\prime }\Vert =\rho _{2}$ and $%
u^{\prime }(t)\geq k_{1}\Vert u^{\prime }\Vert =k_{1}\rho _{2}\geq \rho
,\,t\in \lbrack 0,\,1].$ Similarly, $\Vert v\Vert _{C^{1}}=\Vert v^{\prime
}\Vert =\rho _{2}$ and $v^{\prime }(t)\geq k_{1}\Vert v^{\prime }\Vert
=k_{1}\rho _{2}\geq \rho .$

As in the previous case, by Lemma \ref{Lema k0}, (\ref{eq2}) and (\ref{eq7}) 
\begin{eqnarray*}
T_{1}u(t) &\geq &\int_{\frac{\eta }{\alpha }}^{\eta
}G(t,s)f(s,\,v(s),\,v^{\prime }(s))ds \\
&\geq &k_{0}\int_{\frac{\eta }{\alpha }}^{\eta }g_{0}(s)\text{ }%
f(s,\,v(s),\,v^{\prime }(s))ds\geq k_{0}\xi _{1}\int_{\frac{\eta }{\alpha }%
}^{\eta }g_{0}\text{ }\left( |v(s)|+|v^{\prime }(s)|\right) ds \\
&=&k_{0}\xi _{1}\int_{\frac{\eta }{\alpha }}^{\eta }g_{0}(s)\text{ }ds\int_{%
\frac{\eta }{\alpha }}^{\eta }\left( \left\vert G(t,r)\right\vert
+\left\vert \frac{\partial G}{\partial t}(t,r)\right\vert \right) \left\vert
h(r,\,u(r),\,u^{\prime }(r))\right\vert dr \\
&=&k_{1}k_{0}\Vert u\Vert _{C^{1}}\xi _{1}\xi _{2}\int_{\frac{\eta }{\alpha }%
}^{\eta }g_{0}(s)\text{ }ds\int_{\frac{\eta }{\alpha }}^{\eta }\left(
k_{0}g_{0}(r)+k_{1}g_{1}(r)\right) \text{ }dr>\Vert u\Vert _{C^{1}},
\end{eqnarray*}%
and%
\begin{eqnarray*}
\left( T_{1}u(t)\right) ^{\prime } &\geq &\int_{\frac{\eta }{\alpha }}^{\eta
}\frac{\partial G}{\partial t}(t,s)f(s,\,v(s),\,v^{\prime }(s))ds \\
&\geq &k_{1}\int_{\frac{\eta }{\alpha }}^{\eta }g_{1}(s)\text{ }%
f(s,\,v(s),\,v^{\prime }(s))ds\geq k_{1}\xi _{1}\int_{\frac{\eta }{\alpha }%
}^{\eta }g_{1}(s)\text{ }\left( |v(s)|+|v^{\prime }(s)|\right) ds \\
&\geq &k_{1}\xi _{1}\xi _{2}\int_{\frac{\eta }{\alpha }}^{\eta }g_{1}(s)%
\text{ }ds\int_{\frac{\eta }{\alpha }}^{\eta }\left(
k_{0}g_{0}(r)+k_{1}g_{1}(r)\right) \text{ }k_{1}\Vert u\Vert _{C^{1}}dr \\
&=&\left( k_{1}\right) ^{2}\Vert u\Vert _{C^{1}}\xi _{1}\xi _{2}\int_{\frac{%
\eta }{\alpha }}^{\eta }g_{1}(s)\text{ }ds\int_{\frac{\eta }{\alpha }}^{\eta
}\left( k_{0}g_{0}(r)+k_{1}g_{1}(r)\right) \text{ }dr>\Vert u\Vert _{C^{1}}.
\end{eqnarray*}%
\medskip

\textbf{Case 4.3.} \textit{Suppose that }$\Vert v\Vert \rightarrow +\infty , 
$\textit{\ }$\Vert v^{\prime }\Vert \rightarrow +\infty ,$\ \textit{\ }$%
\Vert u\Vert \rightarrow +\infty $\textit{\ and \ }$\Vert u^{\prime }\Vert
\rightarrow +\infty .$\medskip

Consider $\rho >0$ such that for $(t,v,v^{\prime })\in \lbrack 0,\,1]\times
\lbrack \rho ,\,+\infty )^{2}$ and $(t,\,u,u^{\prime })\in \lbrack
0,\,1]\times \lbrack \rho ,\,+\infty )^{2}$, conditions (\ref{csi1}), (\ref%
{csi2}) and (\ref{eq7}) hold.

Let $u,v\in K$ such that $\Vert u\Vert _{C^{1}}=\rho _{2}$, where $\rho
_{2}:=\max \left\{ 2\rho _{1},\,\frac{\rho }{k_{0}},\,\frac{\rho }{k_{1}}%
\right\} $.

Then $\Vert u\Vert _{C^{1}}=\Vert u\Vert =\Vert u^{\prime }\Vert =\rho _{2}$
and $u(t)\geq k_{0}\Vert u\Vert =k_{0}\rho _{2}\geq \rho ,\,u^{\prime
}(t)\geq k_{1}\Vert u\Vert =k_{1}\rho _{2}\geq \rho ,$ $t\in \lbrack 0,\,1].$
Similarly, $\Vert v\Vert _{C^{1}}=\Vert v\Vert =\Vert v^{\prime }\Vert =\rho
_{2}$ , $v(t)\geq k_{0}\Vert v\Vert =k_{0}\rho _{2}\geq \rho $ and $%
v^{\prime }(t)\geq k_{1}\Vert v^{\prime }\Vert =k_{1}\rho _{2}\geq \rho $ $.$

As before,%
\begin{eqnarray*}
T_{1}u(t) &\geq &k_{0}\xi _{1}\xi _{2}\int_{\frac{\eta }{\alpha }}^{\eta
}g_{0}(s)\text{ }ds\int_{\frac{\eta }{\alpha }}^{\eta }\left(
k_{0}g_{0}(r)+k_{1}g_{1}(r)\right) \text{ }\left( |u(r)|+|u^{\prime
}(r)|\right) dr \\
&\geq &k_{0}\xi _{1}\xi _{2}\int_{\frac{\eta }{\alpha }}^{\eta }g_{0}(s)%
\text{ }ds\int_{\frac{\eta }{\alpha }}^{\eta }\left(
k_{0}g_{0}(r)+k_{1}g_{1}(r)\right) \text{ }\left( k_{0}+k_{1}\right) \Vert
u\Vert _{C^{1}}dr \\
&=&k_{0}\left( k_{0}+k_{1}\right) \Vert u\Vert _{C^{1}}\xi _{1}\xi _{2}\int_{%
\frac{\eta }{\alpha }}^{\eta }g_{0}(s)\text{ }ds\int_{\frac{\eta }{\alpha }%
}^{\eta }\left( k_{0}g_{0}(r)+k_{1}g_{1}(r)\right) \text{ }dr>\Vert u\Vert
_{C^{1}},
\end{eqnarray*}%
and%
\begin{eqnarray*}
\left( T_{1}u(t)\right) ^{\prime } &\geq &k_{1}\int_{\frac{\eta }{\alpha }%
}^{\eta }g_{1}(s)\text{ }f(s,\,v(s),\,v^{\prime }(s))ds\geq k_{1}\xi
_{1}\int_{\frac{\eta }{\alpha }}^{\eta }g_{1}(s)\text{ }\left(
|v(s)|+|v^{\prime }(s)|\right) ds \\
&=&k_{1}\xi _{1}\int_{\frac{\eta }{\alpha }}^{\eta }g_{1}(s)\text{ }ds\int_{%
\frac{\eta }{\alpha }}^{\eta }\left( \left\vert G(t,r)\right\vert
+\left\vert \frac{\partial G}{\partial t}(t,r)\right\vert \right) \left\vert
h(r,\,u(r),\,u^{\prime }(r))\right\vert dr \\
&\geq &k_{1}\xi _{1}\xi _{2}\int_{\frac{\eta }{\alpha }}^{\eta }g_{1}(s)%
\text{ }ds\int_{\frac{\eta }{\alpha }}^{\eta }\left(
k_{0}g_{0}(r)+k_{1}g_{1}(r)\right) \text{ }\left( k_{0}+k_{1}\right) \Vert
u\Vert _{C^{1}}dr \\
&=&k_{1}\left( k_{0}+k_{1}\right) \Vert u\Vert _{C^{1}}\xi _{1}\xi _{2}\int_{%
\frac{\eta }{\alpha }}^{\eta }g_{1}(s)\text{ }ds\int_{\frac{\eta }{\alpha }%
}^{\eta }\left( k_{0}g_{0}(r)+k_{1}g_{1}(r)\right) \text{ }dr>\Vert u\Vert
_{C^{1}}.
\end{eqnarray*}

The other cases follow the same arguments.

Therefore $\Vert T_{1}u\Vert _{C^{1}}\geq \Vert u\Vert _{C^{1}}.$

Then, by Lemma \ref{Lema G-K}, $T_{1}$ has a fixed point in $K\cap (%
\overline{\Omega _{2}}\backslash \Omega _{1}).$

By the same steps it can be proved that $T_{2}$ has a fixed point in $K\cap (%
\overline{\Omega _{2}}\backslash \Omega _{1}),$ too.\medskip

Assume that $(A3)$ and $(A4)$ are verified.\medskip

\textbf{Step 5: \ }$\Vert T_{1}u\Vert _{C^{1}}\geq \Vert u\Vert _{C^{1}}$%
\textit{, for some }$\rho _{3}>0$ \textit{and} \textit{\ }$u\in K\cap
\partial \Omega _{3}$\textit{\ with }$\Omega _{3}=\{u\in E\,:\,\Vert u\Vert
_{C^{1}}<\rho _{3}\}.$\medskip

By $(A3),$ \ it can be chosen $\rho _{3}>0$ such that $(t,\,v,v^{\prime
})\in \lbrack 0,\,1]\times \lbrack 0,\,\rho _{3}]^{2},$ $(t,\,u,u^{\prime
})\in \lbrack 0,\,1]\times \lbrack 0,\,\rho _{3}]^{2}$, and there are $\xi
_{3}$, $\xi _{4}>0$ with%
\begin{eqnarray*}
f(t,v(t),v(t)) &\geq &\xi _{3}\left( |v(t)|+|v^{\prime }(t)||\right) , \\
h(t,u(t),u^{\prime }(t)) &\geq &\xi _{4}\left( |u(t)|+|u^{\prime
}(t)||\right)
\end{eqnarray*}%
and 
\begin{equation}
\min \left\{ 
\begin{array}{c}
\left( k_{0}\right) ^{2}\xi _{3}\xi _{4}\int_{\frac{\eta }{\alpha }}^{\eta
}g_{0}(s)ds\int_{\frac{\eta }{\alpha }}^{\eta }\left(
k_{0}g_{0}(r)+k_{1}g_{1}(r)\right) dr,\text{ } \\ 
k_{0}k_{1}\xi _{3}\xi _{4}\int_{\frac{\eta }{\alpha }}^{\eta }g_{1}(s)\text{ 
}ds\int_{\frac{\eta }{\alpha }}^{\eta }\left(
k_{0}g_{0}(r)+k_{1}g_{1}(r)\right) \text{ }dr \\ 
\xi _{3}\xi _{4}k_{0}k_{1}\int_{\frac{\eta }{\alpha }}^{\eta }g_{0}(s)\text{ 
}ds\int_{\frac{\eta }{\alpha }}^{\eta }\left(
k_{0}g_{0}(r)+k_{1}g_{1}(r)\right) dr \\ 
\left( k_{1}\right) ^{2}\xi _{3}\xi _{4}\int_{\frac{\eta }{\alpha }}^{\eta
}g_{1}(s)\text{ }ds\int_{\frac{\eta }{\alpha }}^{\eta }\left(
k_{0}g_{0}(r)+k_{1}g_{1}(r)\right) \text{ }dr \\ 
\xi _{3}\xi _{4}k_{0}\left( k_{0}+k_{1}\right) \int_{\frac{\eta }{\alpha }%
}^{\eta }g_{0}(s)\text{ }ds\int_{\frac{\eta }{\alpha }}^{\eta }\left(
k_{0}g_{0}(r)+k_{1}g_{1}(r)\right) dr \\ 
k_{1}\left( k_{0}+k_{1}\right) \xi _{3}\xi _{4}\int_{\frac{\eta }{\alpha }%
}^{\eta }g_{1}(s)\text{ }ds\int_{\frac{\eta }{\alpha }}^{\eta }\left(
k_{0}g_{0}(r)+k_{1}g_{1}(r)\right) \text{ }dr%
\end{array}%
\right\} >1.  \label{eq8}
\end{equation}

Let $u\in K$ and $\Vert u\Vert _{C^{1}}=\rho _{3}.$\medskip

\textbf{Case 5.1.} \textit{Suppose }$\Vert u\Vert _{C^{1}}=\Vert u\Vert
=\rho _{3}.$\medskip

By Lemma 2, (\ref{eq2}) and (\ref{eq8}),%
\begin{eqnarray*}
T_{1}u(t) &\geq &\xi _{3}\int_{\frac{\eta }{\alpha }}^{\eta }G(t,s)\left(
|v(s)|+|v^{\prime }(s)|\right) ds \\
&\geq &\xi _{3}\xi _{4}k_{0}\int_{\frac{\eta }{\alpha }}^{\eta }g_{0}(s)%
\text{ }ds\int_{\frac{\eta }{\alpha }}^{\eta }\left(
k_{0}g_{0}(r)+k_{1}g_{1}(r)\right) \text{ }k_{0}\Vert u\Vert _{C^{1}}dr \\
&=&\xi _{3}\xi _{4}\left( k_{0}\right) ^{2}\Vert u\Vert _{C^{1}}\int_{\frac{%
\eta }{\alpha }}^{\eta }g_{0}(s)\text{ }ds\int_{\frac{\eta }{\alpha }}^{\eta
}\left( k_{0}g_{0}(r)+k_{1}g_{1}(r)\right) dr>\Vert u\Vert _{C^{1}},
\end{eqnarray*}%
and%
\begin{eqnarray*}
\left( T_{1}u(t)\right) ^{\prime } &\geq &k_{1}\int_{\frac{\eta }{\alpha }%
}^{\eta }g_{1}(s)\text{ }f(s,\,v(s),\,v^{\prime }(s))ds\geq k_{1}\xi
_{3}\int_{\frac{\eta }{\alpha }}^{\eta }g_{1}(s)\text{ }\left(
|v(s)|+|v^{\prime }(s)|\right) ds \\
&\geq &k_{1}\xi _{3}\xi _{4}\int_{\frac{\eta }{\alpha }}^{\eta }g_{1}(s)%
\text{ }ds\int_{\frac{\eta }{\alpha }}^{\eta }\left(
k_{0}g_{0}(r)+k_{1}g_{1}(r)\right) \text{ }k_{0}\Vert u\Vert _{C^{1}}dr \\
&=&k_{0}k_{1}\Vert u\Vert _{C^{1}}\xi _{3}\xi _{4}\int_{\frac{\eta }{\alpha }%
}^{\eta }g_{1}(s)\text{ }ds\int_{\frac{\eta }{\alpha }}^{\eta }\left(
k_{0}g_{0}(r)+k_{1}g_{1}(r)\right) \text{ }dr>\Vert u\Vert _{C^{1}}.
\end{eqnarray*}%
\medskip

\textbf{Case 5.2.} \textit{Suppose }$\Vert u\Vert _{C^{1}}=\Vert u^{\prime
}\Vert =\rho _{3}.$\medskip

By Lemma 2, (\ref{eq2}) and (\ref{eq8})%
\begin{eqnarray*}
T_{1}u(t) &\geq &\xi _{3}\int_{\frac{\eta }{\alpha }}^{\eta }G(t,s)\left(
|v(s)|+|v^{\prime }(s)|\right) ds \\
&\geq &\xi _{3}k_{0}\int_{\frac{\eta }{\alpha }}^{\eta }g_{0}(s)\int_{\frac{%
\eta }{\alpha }}^{\eta }\left( \left\vert G(t,r)\right\vert +\left\vert 
\frac{\partial G}{\partial t}(t,r)\right\vert \right) \left\vert
h(r,\,u(r),\,u^{\prime }(r))\right\vert drds \\
&\geq &\xi _{3}\xi _{4}k_{0}\int_{\frac{\eta }{\alpha }}^{\eta }g_{0}(s)%
\text{ }ds\int_{\frac{\eta }{\alpha }}^{\eta }\left(
k_{0}g_{0}(r)+k_{1}g_{1}(r)\right) \text{ }k_{1}\Vert u\Vert _{C^{1}}dr \\
&=&\xi _{3}\xi _{4}k_{0}k_{1}\Vert u\Vert _{C^{1}}\int_{\frac{\eta }{\alpha }%
}^{\eta }g_{0}(s)\text{ }ds\int_{\frac{\eta }{\alpha }}^{\eta }\left(
k_{0}g_{0}(r)+k_{1}g_{1}(r)\right) dr>\Vert u\Vert _{C^{1}},
\end{eqnarray*}%
and%
\begin{eqnarray*}
\left( T_{1}u(t)\right) ^{\prime } &\geq &k_{1}\xi _{3}\xi _{4}\int_{\frac{%
\eta }{\alpha }}^{\eta }g_{1}(s)\text{ }ds\int_{\frac{\eta }{\alpha }}^{\eta
}\left( k_{0}g_{0}(r)+k_{1}g_{1}(r)\right) \text{ }\left( |u(r)|+|u^{\prime
}(r)|\right) dr \\
&\geq &k_{1}\xi _{3}\xi _{4}\int_{\frac{\eta }{\alpha }}^{\eta }g_{1}(s)%
\text{ }ds\int_{\frac{\eta }{\alpha }}^{\eta }\left(
k_{0}g_{0}(r)+k_{1}g_{1}(r)\right) \text{ }\left( \underset{r\in \left[ 
\frac{\eta }{\alpha },\,\eta \right] \,}{\min }u(r)+\underset{r\in \left[ 
\frac{\eta }{\alpha },\,\eta \right] \,}{\min }u^{\prime }(r)\right) dr \\
&=&\left( k_{1}\right) ^{2}\Vert u\Vert _{C^{1}}\xi _{3}\xi _{4}\int_{\frac{%
\eta }{\alpha }}^{\eta }g_{1}(s)\text{ }ds\int_{\frac{\eta }{\alpha }}^{\eta
}\left( k_{0}g_{0}(r)+k_{1}g_{1}(r)\right) \text{ }dr>\Vert u\Vert _{C^{1}}.
\end{eqnarray*}

\textbf{Case 5.3.} \textit{Suppose }$\Vert u\Vert _{C^{1}}=\Vert u\Vert
=\Vert u^{\prime }\Vert =\rho _{3}.$\medskip

By Lemma 2, (\ref{eq2}) and (\ref{eq8})

\begin{eqnarray*}
T_{1}u(t) &\geq &\int_{\frac{\eta }{\alpha }}^{\eta
}G(t,s)f(s,\,v(s),\,v^{\prime }(s))ds \\
&\geq &\xi _{3}k_{0}\int_{\frac{\eta }{\alpha }}^{\eta }g_{0}(s)\int_{\frac{%
\eta }{\alpha }}^{\eta }\left( \left\vert G(t,r)\right\vert +\left\vert 
\frac{\partial G}{\partial t}(t,r)\right\vert \right) \left\vert
h(r,\,u(r),\,u^{\prime }(r))\right\vert drds \\
&\geq &\xi _{3}\xi _{4}k_{0}\int_{\frac{\eta }{\alpha }}^{\eta }g_{0}(s)%
\text{ }ds\int_{\frac{\eta }{\alpha }}^{\eta }\left(
k_{0}g_{0}(r)+k_{1}g_{1}(r)\right) \text{ }\left( k_{0}\Vert u\Vert
+k_{1}\Vert u^{\prime }\Vert \right) dr \\
&=&\xi _{3}\xi _{4}k_{0}\left( k_{0}+k_{1}\right) \Vert u\Vert _{C^{1}}\int_{%
\frac{\eta }{\alpha }}^{\eta }g_{0}(s)\text{ }ds\int_{\frac{\eta }{\alpha }%
}^{\eta }\left( k_{0}g_{0}(r)+k_{1}g_{1}(r)\right) dr>\Vert u\Vert _{C^{1}},
\end{eqnarray*}

and%
\begin{eqnarray*}
\left( T_{1}u(t)\right) ^{\prime } &\geq &k_{1}\int_{\frac{\eta }{\alpha }%
}^{\eta }g_{1}(s)\text{ }f(s,\,v(s),\,v^{\prime }(s))ds\geq k_{1}\xi
_{3}\int_{\frac{\eta }{\alpha }}^{\eta }g_{1}(s)\text{ }\left(
|v(s)|+|v^{\prime }(s)|\right) ds \\
&\geq &k_{1}\xi _{3}\xi _{4}\int_{\frac{\eta }{\alpha }}^{\eta }g_{1}(s)%
\text{ }ds\int_{\frac{\eta }{\alpha }}^{\eta }\left(
k_{0}g_{0}(r)+k_{1}g_{1}(r)\right) \text{ }\left( \underset{r\in \left[ 
\frac{\eta }{\alpha },\,\eta \right] \,}{\min }u(r)+\underset{r\in \left[ 
\frac{\eta }{\alpha },\,\eta \right] \,}{\min }u^{\prime }(r)\right) dr \\
&\geq &k_{1}\left( k_{0}+k_{1}\right) \Vert u\Vert _{C^{1}}\xi _{3}\xi
_{4}\int_{\frac{\eta }{\alpha }}^{\eta }g_{1}(s)\text{ }ds\int_{\frac{\eta }{%
\alpha }}^{\eta }\left( k_{0}g_{0}(r)+k_{1}g_{1}(r)\right) \text{ }dr>\Vert
u\Vert _{C^{1}}.
\end{eqnarray*}

In any case, $\Vert T_{1}u\Vert _{C^{1}}\geq \Vert u\Vert _{C^{1}}.\medskip $

\textbf{Step 6: \ }$\Vert T_{1}u\Vert _{C^{1}}\leq \Vert u\Vert _{C^{1}}$%
\textit{, for some }$\rho _{4}>0$ \textit{and} \textit{\ }$u\in K\cap
\partial \Omega _{4}$\textit{\ with }$\Omega _{4}=\{u\in E\,:\,\Vert u\Vert
_{C^{1}}<\rho _{4}\}.\medskip $

Let $u\in K$ and $\Vert u\Vert _{C^{1}}=\rho _{4}.\medskip $

\textbf{Case 6.1.} Suppose that $f$ and$\;h$ are bounded.\medskip

Then there is $N>0$ such that $f(t,v(t),v^{\prime }(t))\leq N$, $%
h(t,u(t),u^{\prime }(t))\leq N,$ $\forall u,\,v\in \lbrack 0,\,\infty ).$

Choose%
\begin{equation*}
\rho _{4}=\max \left\{ 2\rho _{3},\,N\int_{0}^{1}g_{0}(s)ds,\text{ }%
N\int_{0}^{1}g_{1}(s)ds\right\} .
\end{equation*}%
Then%
\begin{equation*}
T_{1}u(t)=\int_{0}^{1}G(t,s)f(s,\,v(s),\,v^{\prime }(s))ds\leq
N\int_{0}^{1}g_{0}(s)ds\leq \rho _{4},\;\text{for }t\in \lbrack 0,\,1],
\end{equation*}%
and%
\begin{equation*}
\left( T_{1}u(t)\right) ^{\prime }=\int_{0}^{1}\frac{\partial G}{\partial t}%
(t,s)f(s,\,v(s),\,v^{\prime }(s))ds\leq N\int_{0}^{1}g_{1}(s)ds\leq \rho
_{4},\;\text{for }t\in \lbrack 0,\,1].
\end{equation*}

Thus, $\Vert T_{1}u\Vert _{C^{1}}\leq \Vert u\Vert _{C^{1}}$. Similarly $%
\Vert T_{2}v\Vert _{C^{1}}\leq \Vert v\Vert _{C^{1}}$ for any $v\in K$ and $%
\Vert v\Vert _{C^{1}}=\rho _{4}$.$\medskip $

\textbf{Case 6.2.} Consider that $f$ is bounded and $h$ is unbounded.$%
\medskip $

So, there is $N>0$ such that $f(t,v(t),v^{\prime }(t))\leq N$, $\forall
\left( v,v^{\prime }\right) \in \lbrack 0,\,+\infty )^{2}$..

By $(A4)$, there exists $M>0$ such that $h(t,u(t),u^{\prime }(t))\leq \mu
\left( |u(t)|+|u^{\prime }(t)|\right) $, whenever $|u(t)|+|u^{\prime }(t)|$ $%
\geq M,$ with $\mu $ verifying

\begin{equation}
\max \left\{ \mu \int_{0}^{1}g_{0}(s)ds,\mu \int_{0}^{1}g_{1}(s)ds\right\} <%
\frac{1}{2}.  \label{miu}
\end{equation}

Setting 
\begin{equation*}
p(r):=\max \{h(t,u(t),u^{\prime }(t))\,:\,t\in \lbrack 0,\,1],\,0\leq u\leq
r,\text{ }0\leq u^{\prime }\leq r\},
\end{equation*}%
we have%
\begin{equation*}
\underset{r\rightarrow \infty }{\lim }p(r)=+\infty .
\end{equation*}

Define%
\begin{equation}
\rho _{4}=\max \left\{ 2\rho _{3},\,M,\,N\int_{0}^{1}g_{0}(s)ds,\text{ }%
N\int_{0}^{1}g_{1}(s)ds\right\} .  \label{eq9}
\end{equation}%
such that $p(\rho _{4})\geq p(r),\;0\leq r\leq \rho _{4}$. Then%
\begin{equation*}
T_{1}u(t)=\int_{0}^{1}G(t,s)f(s,\,v(s),\,v^{\prime }(s))ds\leq
N\int_{0}^{1}g_{0}(s)ds\leq \rho _{4},
\end{equation*}%
and%
\begin{equation*}
\left( T_{1}u(t)\right) ^{\prime }=\int_{0}^{1}\frac{\partial G}{\partial t}%
(t,s)f(s,\,v(s),\,v^{\prime }(s))ds\leq N\int_{0}^{1}g_{1}(s)ds\leq \rho
_{4},\;\text{for }t\in \lbrack 0,\,1].
\end{equation*}

So, $\Vert T_{1}u\Vert _{C^{1}}\leq \Vert u\Vert _{C^{1}}$ for $\Vert u\Vert
_{C^{1}}=\rho _{4}.$

Moreover, if $v\in K$ such that $\Vert v\Vert _{C^{1}}=\rho _{4}$, we have $%
|u(t)|+|u^{\prime }(t)|$ $\geq \rho _{4}\geq M$ , 
\begin{equation}
h(t,u(t),u^{\prime }(t))\leq \mu \left( |u(t)|+|u^{\prime }(t)|\right) \leq
2\mu \rho _{4}  \label{h<2miu}
\end{equation}%
and $p(\rho _{4})\leq 2\mu \rho _{4}.$ Therefore%
\begin{eqnarray*}
T_{2}v(t) &=&\int_{0}^{1}G(t,s)h(s,\,u(s),\,u^{\prime }(s))ds\leq
\int_{0}^{1}g_{0}(s)p(\rho _{4})ds \\
&\leq &p(\rho _{4})\int_{0}^{1}g_{0}(s)ds\leq 2\mu \rho
_{4}\int_{0}^{1}g_{0}(s)ds\leq \rho _{4},
\end{eqnarray*}%
and%
\begin{eqnarray*}
\left( T_{2}v(t)\right) ^{\prime } &=&\int_{0}^{1}\frac{\partial G}{\partial
t}(t,s)h(s,\,u(s),\,u^{\prime }(s))ds\leq \int_{0}^{1}g_{1}(s)p(\rho _{4})ds
\\
&\leq &p(\rho _{4})\int_{0}^{1}g_{1}(s)ds\leq 2\mu \rho
_{4}\int_{0}^{1}g_{1}(s)ds\leq \rho _{4}.
\end{eqnarray*}

Consequently $\Vert T_{2}v\Vert _{C^{1}}\leq \Vert v\Vert _{C^{1}}$ for $%
\Vert v\Vert _{C^{1}}=\rho _{4}.\medskip $

\textbf{Case 6.3.} Suppose that $f$ is unbounded and $h$ is bounded.$%
\medskip $

Then, there is $N>0$ such that $h(t,u(t),u^{\prime }(t))\leq N$, $\forall
\left( u,u^{\prime }\right) \in \lbrack 0,\,+\infty )^{2}$., and, by $(A4)$,
there exists $M>0$ such that $f(t,v(t),v^{\prime }(t))\leq \mu \left(
|v|+|v^{\prime }|\right) $, for $\left( |v|+|v^{\prime }|\right) \geq M$
with $\mu $ satisfying (\ref{miu}).

Choosing $\rho _{4}$ as in (\ref{eq9}), the arguments follow like in the
previous case.$\medskip $

\textbf{Case 6.4.} Consider that $f$ and $h$ are unbounded.$\medskip $

By $(A4)$, there is $M>0$ such that $f(t,v(t),v^{\prime }(t))\leq \mu \left(
|v|+|v^{\prime }|\right) $, $h(t,u(t),u^{\prime }(t))\leq \mu \left(
|u|+|u^{\prime }||\right) $ for $|v|+|v^{\prime }|\geq M$ and $%
|u|+|u^{\prime }|\geq $ $M$ with $\mu $ as in (\ref{miu}).

Setting 
\begin{eqnarray*}
p(r) &:&=\max \{h(t,u(t),u^{\prime }(t))\,:\,t\in \lbrack 0,\,1],\,0\leq
u\leq r,\text{ }0\leq u^{\prime }\leq r\}, \\
q(r) &:&=\max \{f(t,v(t),v^{\prime }(t)):\,t\in \lbrack 0,\,1],\,0\leq v\leq
r,\text{ }0\leq v^{\prime }\leq r\}
\end{eqnarray*}%
we have%
\begin{equation*}
\underset{r\rightarrow \infty }{\lim }p(r)=+\infty \text{ and }\underset{%
r\rightarrow \infty }{\lim }q(r)=+\infty .
\end{equation*}%
Choose 
\begin{equation*}
\rho _{4}=\max \left\{ 2\rho _{3},\,M\right\}
\end{equation*}%
such that $p(\rho _{4})\geq p(r)$ and $q(\rho _{4})\geq q(r)$ for$\;0\leq
r\leq \rho _{4}$.

Let $u,\,v\in K$ and $\Vert u\Vert _{C^{1}}=\Vert v\Vert _{C^{1}}=\rho _{4}$.

Arguing as in (\ref{h<2miu}) it can be easily shown that $\Vert T_{1}u\Vert
_{C^{1}}\leq \Vert u\Vert _{C^{1}},\,\Vert T_{2}v\Vert _{C^{1}}\leq \Vert
v\Vert _{C^{1}}$. $\medskip $

By Lemma \ref{Lema G-K} the operators $T_{1},\;T_{2}$ has\ a fixed point in $%
K\cap (\overline{\Omega _{4}}\backslash \Omega _{3}),$ therefore $T=\left(
T_{1},T_{2}\right) $ has a fixed point $(u,v)$ which is a positive solution
of the initial problem.

Moreover these functions $u$ and $v$ are given by 
\begin{equation*}
\left\{ 
\begin{array}{c}
u(t)=\int_{0}^{1}G(t,s)f(s,\,v(s),\,v^{\prime }(s))ds \\ 
\\ 
v(t)=\int_{0}^{1}G(t,s)h(s,\,u(s),\,u^{\prime }(s))ds.%
\end{array}%
\right.
\end{equation*}%
and are both increasing functions.
\end{proof}

\section{Example}

Consider the following third order nonlinear system%
\begin{equation}
\left\{ 
\begin{array}{c}
-u^{\prime \prime \prime }(t)=\left( t^{2}+1\right) \left( e^{-v(t)}+\sqrt{%
\left\vert v^{\prime }(t)\right\vert }\right) \\ 
\\ 
-v^{\prime \prime \prime }(t)=\left( \,u(t)+1\right) ^{2}\arctan \left(
\left\vert u^{\prime }(t)\right\vert +1\right) \\ 
\\ 
u(0)=u^{\prime }(0)=0,u^{\prime }(1)=\frac{3}{2}u^{\prime }\left( \frac{1}{2}%
\right) \\ 
\\ 
v(0)=v^{\prime }(0)=0,v^{\prime }(1)=\frac{3}{2}v^{\prime }\left( \frac{1}{2}%
\right) .%
\end{array}%
\right.  \label{EX}
\end{equation}

In fact this problem is a particular case of system (\ref{eq1}) with%
\begin{eqnarray*}
f(t,\,v(t),\,v^{\prime }(t)) &:&=\left( t^{2}+1\right) \left( e^{-v(t)}+%
\sqrt{\left\vert v^{\prime }(t)\right\vert }\right) \\
h(t,\,u(t),\,u^{\prime }(t)) &:&=\left( \,u(t)+1\right) ^{2}\arctan \left(
\left\vert u^{\prime }(t)\right\vert +1\right) , \\
\eta &=&\frac{1}{2}\text{ and \ }\alpha =\frac{3}{2}.
\end{eqnarray*}

It can be easily check that the above functions are non-negative and verify
the assumptions $(A3)$ and $(A4).$

$f$ and $h$ are non-negative functions, because they are product of non
negative functions. Note that, $\forall t$ and $\forall (u(t),\;v(t))\in %
\big(C^{3}[0,\,1],\,(0,\;+\infty )$, $(t^2+1)\geq1$, $e^{-v(t)}=\frac{1}{%
e^{v(t)}}\geq0$, $\sqrt{\left\vert v^{\prime }(t)\right\vert }\geq0$, $%
\left( \,u(t)+1\right) ^{2}\geq1$ and for definition of arc tangent
function, we know $\arctan:\;\mathbb{R}\rightarrow \left]-\frac{\pi}{2},\;%
\frac{\pi}{2}\right[$ such that $\arctan(x)\rightarrow\frac{\pi}{2}$, as $%
x\rightarrow+\infty$ and $\arctan(x)\rightarrow-\frac{\pi}{2}$, as $%
x\rightarrow-\infty$. So, $f\geq0$, $h\geq0$.

Finally, as

\begin{eqnarray*}
\underset{t\in \lbrack 0,1],\text{ }\Vert v\Vert _{C^{1}}\rightarrow 0}{%
\liminf }\frac{\left( t^{2}+1\right) \left( e^{-v(t)}+\sqrt{\left\vert
v^{\prime }(t)\right\vert }\right) }{|v|+|v^{\prime }|} &=&+\infty , \\
\underset{t\in \lbrack 0,1],\text{ }\Vert u\Vert _{C^{1}}\rightarrow 0}{%
\liminf }\frac{\left( \,u(t)+1\right) ^{2}\arctan \left( \left\vert
u^{\prime }(t)\right\vert +1\right) }{|u|+|u^{\prime }|} &=&+\infty
\end{eqnarray*}%
conditionn $(A3)$ kolds and

\begin{eqnarray*}
\underset{t\in \lbrack 0,1],\text{ }\Vert v\Vert _{_{C^{1}}}\rightarrow
+\infty }{\limsup }\frac{\left( t^{2}+1\right) \left( e^{-v(t)}+\sqrt{%
\left\vert v^{\prime }(t)\right\vert }\right) }{|v|+|v^{\prime }|} &=&0, \\
\underset{t\in \lbrack 0,1],\text{ }\Vert u\Vert _{_{C^{1}}}\rightarrow
+\infty }{\limsup }\frac{\left( \,u(t)+1\right) ^{2}\arctan \left(
\left\vert u^{\prime }(t)\right\vert +1\right) }{|u|+|u^{\prime }|} &=&0
\end{eqnarray*}%
assumptiom $(A4)$ is satisfied.

Therefore, by Theorem \ref{Main thm}, problem (\ref{EX}) has at least a
positive solution $(u(t),\,v(t))\in \left( C^{3}[0,\,1]\right) ^{2}$, that
is $u(t)>0$, $v(t)>0,$ $\forall t\in \lbrack 0,\,1]$.

\end{document}